\documentclass[12pt,reqno]{amsart}

\usepackage{amsmath,amssymb,amsthm,mathtools}
\usepackage[margin=0.91in]{geometry}
\usepackage{enumitem}
\usepackage[hidelinks]{hyperref}
\usepackage{microtype}
\usepackage{needspace}

\numberwithin{equation}{section}

\newtheorem{theorem}{Theorem}[section]
\newtheorem{lemma}[theorem]{Lemma}
\newtheorem{proposition}[theorem]{Proposition}

\theoremstyle{definition}
\newtheorem{definition}[theorem]{Definition}
\theoremstyle{remark}
\newtheorem{remark}[theorem]{Remark}

\newcommand{\PSL}{\operatorname{PSL}}
\newcommand{\PGL}{\operatorname{PGL}}
\newcommand{\SL}{\operatorname{SL}}
\newcommand{\GL}{\operatorname{GL}}
\newcommand{\PSU}{\operatorname{PSU}}
\newcommand{\PSp}{\operatorname{PSp}}
\newcommand{\PSO}{\operatorname{PSO}}
\newcommand{\POmega}{\operatorname{P\Omega}}
\newcommand{\F}{\mathbb F}
\newcommand{\Dcal}{\mathcal D}
\newcommand{\Ecal}{\mathcal E}

\title[Rohrbach's conjecture]{Rohrbach's Conjecture for Finite Groups}
\author{Yong Yang}
\address{Department of Mathematics, Texas State University, 601 University Drive, San Marcos, TX 78666, USA}
\email{yang@txstate.edu}
\date{}


\begin{document}

\begin{abstract}
Rohrbach's conjecture holds for finite groups. In fact, there is an absolute constant $C$ such that, for every finite group $G$, every positive integer $h$, and all nonnegative real numbers $\alpha_1,\ldots,\alpha_h$ with $\alpha_1+\cdots+\alpha_h=1$, there are subsets $A_1,\ldots,A_h\subseteq G$ satisfying
\[
G=A_1\cdots A_h,
\qquad
|A_1|\leq |G|^{\alpha_1},
\qquad
|A_i|\leq C|G|^{\alpha_i}\quad(2\leq i\leq h).
\]
Consequently, for every $h$ there is a constant $c_h$ such that every finite group $G$ has a subset $A$ with
\[
A^h=G\qquad\text{and}\qquad |A|\leq c_h|G|^{1/h}.
\]
\end{abstract}

\maketitle

\section{Introduction}

Let $G$ be a finite group and let $h\geq2$. A subset $A\subseteq G$ is an $h$-basis if $A^h=G$. Since $|G|\leq |A|^h$, every $h$-basis satisfies $|A|\geq |G|^{1/h}$. Rohrbach asked in 1937 whether, for each fixed $h$, there is a constant $c_h$ such that every finite group $G$ has an $h$-basis of size at most $c_h|G|^{1/h}$ \cite{Rohrbach1,Rohrbach2}.

The case $h=2$ was proved for all finite groups in \cite[Theorems~1 and~2]{KozmaLev1992}. For arbitrary $h$, Jia proved Rohrbach's conjecture for finite nilpotent groups \cite{Jia}. Kozma and Lev subsequently proved a stronger weighted result for finite solvable and alternating groups \cite{KozmaLev}. The result below concerns arbitrary finite groups, with a weighted constant independent of both $G$ and $h$.

More generally, \cite{KozmaLev} proves that if every composition factor of $G$ is cyclic or alternating, then for every choice of nonnegative $\alpha_1,\ldots,\alpha_h$ with sum $1$ there are subsets $A_1,\ldots,A_h\subseteq G$ such that
\[
G=A_1\cdots A_h,\qquad
|A_1|\leq |G|^{\alpha_1},\qquad
|A_i|\leq2|G|^{\alpha_i}\quad(2\leq i\leq h).
\]
The extension lemma in that paper reduces the uniform decomposition problem for finite groups to finite simple groups \cite[p.~390]{KozmaLev}. Section~2 gives the same reduction with an arbitrary fixed constant.

Babai, G\'al, Kimmel and Lokam later reported an extension of Pyber covering groups with linear and symplectic composition factors, with a suitable absolute constant in place of $2$ \cite[Section~7]{BabaiGalKimmelLokam}. They cite this result as a private communication. Their paper formulated the Modified Rohrbach Problem, asking for one absolute constant in uniform $h$-decompositions of arbitrary finite groups.

The main result completes this program, and in fact gives the stronger weighted conclusion above.

\begin{theorem}\label{thm:global}
There is an absolute constant $C\geq1$ such that the following holds. Let $G$ be a finite group, let $h\geq1$ be an integer, and let $\alpha_1,\ldots,\alpha_h\geq0$ satisfy $\sum_i\alpha_i=1$. Then there are subsets $A_1,\ldots,A_h\subseteq G$ such that
\[
G=A_1\cdots A_h,
\]
\[
|A_1|\leq |G|^{\alpha_1},
\qquad
|A_i|\leq C|G|^{\alpha_i}\quad(2\leq i\leq h).
\]
\end{theorem}

\begin{theorem}\label{thm:rohrbach}
Rohrbach's conjecture holds for all finite groups. More precisely, for every positive integer $h$ there is a constant $c_h$ such that every finite group $G$ has a subset $A\subseteq G$ satisfying
\[
A^h=G,
\qquad
|A|\leq c_h|G|^{1/h}.
\]
One may take $c_h=1+C(h-1)$, where $C$ is the constant in Theorem~\ref{thm:global}.
\end{theorem}

The main step is the following uniform theorem for groups of Lie type.

\begin{theorem}\label{thm:lie}
There is an absolute constant $C_{\rm Lie}$ such that every nonabelian finite simple group of Lie type satisfies the conclusion of Theorem~\ref{thm:global} with $C=C_{\rm Lie}$.
\end{theorem}

Section~2 contains the decomposition lemmas used throughout. Section~3 treats the linear groups by a Singer-cycle transversal with an explicit determinant correction. Section~4 establishes the Bruhat transversal used in the remaining unbounded-rank classical families. Sections~5 and~6 treat the unitary, symplectic, and orthogonal groups case by case. Section~7 handles the exceptional groups. The proof of the main theorem is completed in Section~8.

\section{General decomposition lemmas}

For finite subsets $Y_1,\ldots,Y_s$ of a group, the product $Y_1\cdots Y_s$ is called \emph{exact} if the multiplication map
\[
Y_1\times\cdots\times Y_s\longrightarrow Y_1\cdots Y_s
\]
is bijective.

The following covering form of the weighted decomposition property is convenient.

\begin{definition}\label{def:C}
Let $C\geq1$, and let $X$ be a nonempty finite subset of a group $\Gamma$. The set $X$ has property $\Dcal_C$ if, whenever $r\geq1$ is an integer and $x_1,\ldots,x_r\geq1$ satisfy
\[
x_1\cdots x_r=|X|,
\]
there are finite subsets $X_1,\ldots,X_r\subseteq\Gamma$ such that
\[
X\subseteq X_1\cdots X_r,
\]
\[
|X_1|\leq x_1,
\qquad
|X_i|\leq Cx_i\quad(2\leq i\leq r).
\]
\end{definition}

For a finite group $H$, property $\Dcal_C$ always refers to coverings by subsets of $H$ itself. If $|H|>1$, the substitutions
\[
x_i=|H|^{\alpha_i},\qquad
\alpha_i=\frac{\log x_i}{\log|H|}
\]
identify the group version of Definition~\ref{def:C} with the weighted statement of Theorem~\ref{thm:global}. The trivial group satisfies both formulations directly. For a subset of a larger group, the ambient group is specified in Definition~\ref{def:C}. The covering property is invariant under left and right translations, and under replacing $C$ by a larger constant.

\begin{lemma}\label{lem:concat}
Let $C\geq1$, and let $X$ and $Y$ be finite subsets of the same group, each with property $\Dcal_C$. Set $M=|X||Y|$. If $x_1,\ldots,x_r\geq1$ and $x_1\cdots x_r=M$, then there are subsets $A_1,\ldots,A_r$ such that
\[
XY\subseteq A_1\cdots A_r,
\]
\[
|A_1|\leq x_1,
\qquad
|A_i|\leq Cx_i\quad(2\leq i\leq r).
\]
The same assertion holds for any finite product $Y_1\cdots Y_s$, with $M=\prod_j|Y_j|$.
\end{lemma}

\begin{proof}
Choose $j$ so that
\[
x_1\cdots x_{j-1}\leq |X|\leq x_1\cdots x_j.
\]
Write $x_j=ab$, where
\[
x_1\cdots x_{j-1}a=|X|,
\qquad
b x_{j+1}\cdots x_r=|Y|.
\]
Apply $\Dcal_C$ to $X$ with the numbers $x_1,\ldots,x_{j-1},a$, and to $Y$ with the numbers $b,x_{j+1},\ldots,x_r$. At the joining position multiply the last subset for $X$ by the first subset for $Y$. The first subset in the decomposition of $Y$ has no factor $C$, so the joining position still has coefficient at most $C$. If $j=1$, both joining subsets have coefficient $1$.

For the assertion about several factors, subdivide an interval of length $\log M$ first into consecutive intervals of lengths $\log|Y_j|$, and then into consecutive intervals of lengths $\log x_i$. Each interval belonging to $Y_j$ is subdivided at its intersections with the second partition. Apply $\Dcal_C$ to $Y_j$ with the exponentials of these subinterval lengths, and multiply the resulting subsets in their original order within each interval belonging to $x_i$. A whole $Y_j$ needs no enlargement. An interval belonging to $x_i$ contains at most one part that is not the first part of its $Y_j$; only that part can contribute a factor $C$. The first interval contains no such part. This proves the claimed bounds even when the original product is not exact. Singleton $Y_j$ may be incorporated into a neighbouring factor by translation, and a prescribed value $x_i=1$ may be represented by $\{1\}$. If $M=1$, the assertion is immediate. If every covering used for the individual factors contains $1$ in each position, the resulting covering has the same property.
\end{proof}

A \emph{cyclic progression} is a set $\{1,x,\ldots,x^{m-1}\}$ in which the displayed powers are distinct. The following elementary version of the cyclic covering argument suffices here.

\begin{lemma}\label{lem:progression}
Every cyclic progression has property $\Dcal_2$. The covering subsets may all be chosen to contain $1$.
\end{lemma}

\begin{proof}
Let $P=\{1,x,\ldots,x^{m-1}\}$ and let $x_1\cdots x_r=m$, where every $x_i\geq1$. Define integers
\[
b_1=\lfloor x_1\rfloor,\qquad
b_i=\lfloor2x_i\rfloor\quad(2\leq i\leq r).
\]
For real numbers $a,b\geq1$,
\[
\frac{(a+1)(b+1)}2\leq ab+1,
\]
since the difference between the right and left sides is $(a-1)(b-1)/2$. Repeated application gives
\[
m<\frac{\prod_{i=1}^r(b_i+1)}{2^{r-1}}
  \leq \prod_{i=1}^r b_i+1.
\]
Both $m$ and $\prod_i b_i$ are integers, so $\prod_i b_i\geq m$.

Define $a_1=1$ and $a_i=b_1\cdots b_{i-1}$ for $i>1$, and let
\[
P_i=\{x^{j a_i}:0\leq j<b_i\}.
\]
The corresponding sums of exponents contain every integer from $0$ to $\prod_i b_i-1$. Hence $P\subseteq P_1\cdots P_r$. Moreover, $1\in P_i$ and $|P_i|\leq b_i$, which are precisely the required bounds.
\end{proof}

The next lemma extends the solvable case to a subgroup together with one additional element.

\begin{lemma}\label{lem:solvable-plus-one}
Let $U$ be a finite solvable subgroup of a group and let $z\notin U$. Then
\[
U\cup\{z\}
\]
has property $\Dcal_3$.
\end{lemma}

\begin{proof}
The assertion is immediate when $r=1$. Assume $r\geq2$, and let $x_1\cdots x_r=|U|+1$. Choose numbers $y_i$ with
\[
1\leq y_i\leq x_i,
\qquad
\prod_i y_i=|U|.
\]
If $|U|=1$, take every $y_i=1$. Otherwise set
\[
t=\frac{\log |U|}{\log(|U|+1)}
\qquad\text{and}\qquad
y_i=x_i^t.
\]
Then $1\leq y_i\leq x_i$ and $\prod_i y_i=|U|$. A composition series of $U$ gives an exact product of cyclic progressions, each containing $1$: at a composition factor of prime order $p_0$, lift a generator and take its first $p_0$ powers as a transversal. Lemma~\ref{lem:progression} and Lemma~\ref{lem:concat} therefore give subsets $B_i$ containing $1$ such that
\[
U\subseteq B_1\cdots B_r,
\qquad 1\in B_i,
\]
with
\[
|B_1|\leq y_1,
\qquad
|B_i|\leq2y_i\quad(i\geq2).
\]
Replace $B_2$ by $B_2\cup\{z\}$. The product now contains $U\cup\{z\}$ and
\[
|B_2\cup\{z\}|\leq2y_2+1\leq3x_2.
\]
All other bounds are unchanged.
\end{proof}

\begin{lemma}\label{lem:normal-extension}
Let $C\geq1$ and let $N\lhd G$. If both $N$ and $G/N$ have property $\Dcal_C$, then $G$ has property $\Dcal_C$.
\end{lemma}

\begin{proof}
Let $x_1,\ldots,x_r\geq1$ with $x_1\cdots x_r=|G|$, and write $n=|N|$. Choose $j$ so that
\[
x_1\cdots x_{j-1}\leq n\leq x_1\cdots x_j.
\]
Write $x_j=ab$, where
\[
x_1\cdots x_{j-1}a=n,
\qquad
b x_{j+1}\cdots x_r=|G/N|.
\]
Apply $\Dcal_C$ to $N$ with the first list and to $G/N$ with the second. Choose representatives in $G$ for the subsets obtained in the quotient. If $R_j,\ldots,R_r$ are these lifted subsets, then $N R_j\cdots R_r=G$: their product meets every coset of $N$, and multiplication by $N$ gives the full preimage. Multiplying the last subset for $N$ by $R_j$ therefore gives a covering of $G$. At the joining position the quotient subset is the first subset in its decomposition and therefore carries no factor $C$. Hence the joining position has coefficient at most $C$, and all remaining bounds are unchanged. If $j=1$, neither of the two joining subsets carries a factor $C$. Thus $G$ has property $\Dcal_C$.
\end{proof}

This is the extension argument in \cite[Lemma~2]{KozmaLev}. In particular, once one absolute constant works for all finite simple groups, the same constant works for all finite groups by induction along a composition series.

\begin{lemma}\label{lem:index}
Let $C\geq1$, and let $M\leq G$ have index $d$. If $G$ has property $\Dcal_C$, then $M$ has property $\Dcal_{dC}$.
\end{lemma}

\begin{proof}
The case of one factor is immediate. Let $r\geq2$ and let $x_1\cdots x_r=|M|$. Apply $\Dcal_C$ in $G$ to
\[
x_1,x_2,\ldots,x_{r-1},d x_r.
\]
There are subsets $A_i\subseteq G$ with $G=A_1\cdots A_r$ and
\[
|A_1|\leq x_1,\qquad
|A_i|\leq Cx_i\quad(2\leq i<r),\qquad
|A_r|\leq Cd x_r.
\]
Choose a transversal $R$ for the right cosets of $M$ in $G$, containing $1$. For $g\in G$, let $\overline g\in R$ represent $Mg$. Define
\[
B_1=\{a\,\overline a^{-1}:a\in A_1\},
\]
\[
B_i=\{ua\,\overline{ua}^{-1}:u\in R,\ a\in A_i\}
\qquad(2\leq i<r),
\]
and
\[
B_r=\{\overline{a^{-1}}a:a\in A_r\}.
\]
Every displayed element belongs to $M$. For the last set, this follows from
$M\overline{a^{-1}}=Ma^{-1}$. Hence
\[
|B_1|\leq x_1,\qquad
|B_i|\leq dC x_i\quad(2\leq i<r),\qquad
|B_r|\leq |A_r|\leq dC x_r.
\]

For a word $a_1\cdots a_r\in M$ with $a_i\in A_i$, define $u_0=1$ and $u_i=\overline{u_{i-1}a_i}$. Inductively, $Mu_i=Ma_1\cdots a_i$, so $u_r=1$. Cancellation gives
\[
a_1\cdots a_r
=\prod_{i=1}^r\bigl(u_{i-1}a_i u_i^{-1}\bigr).
\]
The first factor lies in $B_1$, and each intermediate factor lies in the corresponding $B_i$. Since $u_{r-1}a_r\in M$, the unique possible representative $u_{r-1}$ is $\overline{a_r^{-1}}$. Thus the last factor belongs to $B_r$. It follows that $M=B_1\cdots B_r$, proving the lemma. No normality of $M$ is required.
\end{proof}

A finite subset $X$ of a group has property $\Ecal$ if $X$ is an exact product of cyclic progressions.

\begin{lemma}\label{lem:E-extension}
For finite groups, property $\Ecal$ is preserved under extensions.
\end{lemma}

\begin{proof}
Let $N\lhd G$, and assume that $N$ and $G/N$ have property $\Ecal$. Lift a generator from each cyclic progression in an exact decomposition of $G/N$. The lifted progressions give one representative for every element of $G/N$. Appending an exact cyclic-progression decomposition of $N$ gives a unique expression for every element of $G$.
\end{proof}

Every finite solvable group has property $\Ecal$, by Lemma~\ref{lem:E-extension} applied along a composition series with prime-order factors.

\begin{lemma}\label{lem:solvable-transversal}
Let $J\leq H$, where $H$ is finite and solvable. There is a left transversal $T$ for $J$ in $H$ which has property $\Ecal$.
\end{lemma}

\begin{proof}
The proof is by induction on $|H|$. If $H=J$, the transversal $\{1\}$ has the required property. Otherwise $H\neq1$; let $N$ be a minimal normal subgroup of $H$. Then $N$ is elementary abelian. By induction in $H/N$, a transversal for $JN/N$ in $H/N$ is an exact product of cyclic progressions. Lifting the generators of those progressions gives a set $T_0\subseteq H$, itself an exact product of cyclic progressions, which is a transversal for $JN$ in $H$.

The subgroup $N\cap J$ is a vector subspace of the elementary abelian group $N$. Choose a vector-space complement $V$ to $N\cap J$ in $N$. Then
\[
N=V(N\cap J),\qquad V\cap J=1.
\]
Thus $JN=VJ$, and multiplication $V\times J\longrightarrow JN$ is bijective. A basis of $V$ writes $V$ as an exact product of cyclic subgroups of prime order. Consequently
\[
T=T_0V
\]
is a left transversal for $J$ in $H$ and is an exact product of cyclic progressions.
\end{proof}

Lemma~\ref{lem:progression} and Lemma~\ref{lem:concat} show that every set with property $\Ecal$ has property $\Dcal_2$. Since the progression coverings in Lemma~\ref{lem:progression} contain $1$, the same is true for the resulting coverings of a set with property $\Ecal$.

\section{Linear groups}

For the linear family, a small determinant correction to the usual Singer cycle avoids passage to $\PGL_n(q)$.

\begin{proposition}\label{prop:linear}
For every $n\geq2$ and every prime power $q$, the group $\PSL_n(q)$ has property $\Dcal_2$.
\end{proposition}

\begin{proof}
Let
\[
G=\PSL_n(q),\qquad \widehat G=\PGL_n(q),\qquad
d=(n,q-1),\qquad v=\frac{q^n-1}{q-1}.
\]
Actions are on the left, with $(gh)\xi=g(h\xi)$ for group elements $g,h$ and a vector or point $\xi$. Identify the natural module with $\F_{q^n}$, and let $\alpha$ generate $\F_{q^n}^{\times}$. Multiplication by $\alpha$ induces a Singer element $c\in\widehat G$ of order $v$, and $\langle c\rangle$ acts regularly on the projective points. The determinant of multiplication by $\alpha$ is its norm to $\F_q$, which generates $\F_q^{\times}$. Since
\[
v=1+q+\cdots+q^{n-1}\equiv n\pmod{q-1},
\]
one has $d\mid v$. Hence the image of $c$ in
\[
\widehat G/G\cong \F_q^{\times}/(\F_q^{\times})^n
\]
has order $d$. Therefore
\[
\Sigma=\langle c\rangle\cap G=\langle c^d\rangle
\]
has order $v/d$.

Let $\omega=\F_q1$ be the projective point represented by $1\in\F_{q^n}$. Since $d\leq n$, the vectors
\[
1,\alpha,\ldots,\alpha^{d-1}
\]
are linearly independent over $\F_q$: the element $\alpha$ has degree $n$, since it generates $\F_{q^n}^{\times}$. An element $A\in\SL_n(q)$ may be chosen so that its projective image $a\in G$ satisfies
\[
a^r\omega=c^r\omega\qquad(0\leq r<d).
\]
If $d<n$, extend the displayed vectors to a basis, let $A$ cyclically permute the first $d$ basis vectors, multiply one remaining basis vector by $(-1)^{d-1}$, and fix the other basis vectors. The determinant is then $1$. If $d=n$, define
\[
A(\alpha^i)=\alpha^{i+1}\quad(0\leq i<n-1),\qquad
A(\alpha^{n-1})=(-1)^{n-1}.
\]
Again $\det A=1$, and the asserted projective equalities hold. In both constructions $A^d$ is scalar, so $a^d=1$.

Let
\[
T=\{1,a,\ldots,a^{d-1}\}.
\]
For $0\leq k<v/d$ and $0\leq r<d$, the element $(c^d)^ka^r\in G$ sends $\omega$ to
\[
c^{dk+r}\omega.
\]
These are all $v$ projective points, each exactly once. In particular, $T=\langle a\rangle$ is cyclic of order $d$. Let $P$ be the stabilizer of $\omega$ in $G$. Multiplication therefore gives an exact product
\[
G=\Sigma TP.
\]
Both $\Sigma$ and $T$ are cyclic progressions and have property $\Dcal_2$.

Induction on $n$ completes the proof. For $n=2$, the point stabilizer $P$ is solvable, including when $q=2$ or $3$. Assume $n\geq3$. Let $Q$ be the unipotent radical of $P$ and let $L=P/Q$. Before passage to the projective group, the Levi subgroup of the line stabilizer consists of
\[
\operatorname{diag}\bigl((\det D)^{-1},D\bigr),
\qquad D\in\GL_{n-1}(q).
\]
Consequently
\[
L\cong\GL_{n-1}(q)/Z,\qquad
Z=\{\lambda I_{n-1}:\lambda\in\F_q^{\times},\ \lambda^d=1\}.
\]
The group $Q$ has order $q^{n-1}$. The image $L_0$ of $\SL_{n-1}(q)$ is normal in $L$, the quotient $L/L_0$ is cyclic, and $L_0$ is a central quotient of $\SL_{n-1}(q)$. Hence $L_0$ has a central cyclic subgroup $Z_0$ with
\[
L_0/Z_0\cong\PSL_{n-1}(q).
\]
The induction hypothesis, the solvable case, and Lemma~\ref{lem:normal-extension} therefore give property $\Dcal_2$ first for $L_0$, then for $L$, and finally for $P$. Since $\Sigma TP$ is exact, Lemma~\ref{lem:concat} gives property $\Dcal_2$ for $G$.
\end{proof}

\begin{remark}\label{rem:linear-levi}
Let $m\geq2$, let $q_0$ be a prime power, and let $L$ be a central quotient of a subgroup $H$ satisfying $\SL_m(q_0)\leq H\leq\GL_m(q_0)$. Let $L_0$ be the image of $\SL_m(q_0)$ in $L$. Then $L_0\lhd L$, the quotient $L/L_0$ is cyclic, and $L_0$ has a central cyclic subgroup $Z_0$ such that
\[
L_0/Z_0\cong\PSL_m(q_0).
\]
Proposition~\ref{prop:linear}, the solvable case, and Lemma~\ref{lem:normal-extension} give property $\Dcal_2$ for $L$. This also applies to the linear sections of the classical Levi quotients below.
\end{remark}

\section{Bruhat transversals}

Let $G$ be a finite group of Lie type with its standard $BN$-pair $(B,\mathcal N)$, where $B=UT$, $U$ is the unipotent radical and $T=B\cap\mathcal N$ is the torus in the $BN$-pair. Let $U^-$ be the opposite unipotent subgroup. The relative Weyl group is
\[
W=\mathcal N/T,
\]
with Coxeter length $\ell$. Here $U$ is a $p$-group, where $p$ is the defining characteristic, and $T$ is abelian; hence $B$ is solvable. For each $w\in W$, choose a representative $\dot w\in\mathcal N$ and define
\[
U_w=U\cap\dot wU^-\dot w^{-1}.
\]
The torus $T$ normalizes $U$ and $U^-$, and therefore normalizes each $U_w$. In particular, the subgroup $U_w$ is independent of the chosen representative. These are subgroups of prime-power order. The refined Bruhat decomposition gives a bijection
\[
U_w\times B\longrightarrow B\dot wB,
\qquad (u,b)\longmapsto u\dot wb.
\]
If $\ell(xy)=\ell(x)+\ell(y)$, then
\[
U_{xy}=U_x(\dot xU_y\dot x^{-1})
\]
is an exact product. The representative $\dot x\dot y$ differs from $\dot{xy}$ by an element of $T$. Since $T$ normalizes every unipotent parameter subgroup, such a torus factor may be moved through subsequent factors and absorbed in the right-hand $B$-factor. See \cite[Chapter~2]{Carter} and \cite[Chapters~1 and~3]{DigneMichel}.

These are statements for the standard finite $BN$-pair itself. In particular, they apply without change to the twisted groups when $W$ is the relative Weyl group and $U_w$ is the corresponding relative unipotent parameter subgroup. Passing to a central quotient does not alter the bijections. Indeed, let $Z\leq T$ be the central kernel and let $P$ be a standard parabolic subgroup. If $u\dot wP$ and $u'\dot wP$ have the same image in the quotient, then $u\dot wPZ=u'\dot wPZ$. Since $Z\leq P$, these are already equal cosets before taking the quotient. Moreover, $Z$ intersects every unipotent parameter subgroup trivially. The same argument applies with $P=B$.

For a standard parabolic $P$ with Weyl subgroup $W_J$, let $W^J$ be the set of minimal-length representatives of the cosets $wW_J$. Then
\[
G=\bigsqcup_{w\in W^J}U_w\dot wP,
\]
and the multiplication map $U_w\times P\longrightarrow U_w\dot wP$ is bijective for every $w\in W^J$. Thus different $w$ describe disjoint $B$-orbits on $G/P$, while the elements of $U_w$ distinguish the individual cosets in the orbit.

\Needspace{11\baselineskip}
\begin{lemma}\label{lem:Bweyl}
Let $W=W(B_m)$ with simple reflections $s_1,\ldots,s_m$, where $\langle s_1,\ldots,s_{m-1}\rangle\cong S_m$. Define
\[
v_i=s_{m-i+1}s_{m-i+2}\cdots s_m
\qquad(1\leq i\leq m).
\]
Then the $2^m$ elements
\[
v_1^{\varepsilon_1}\cdots v_m^{\varepsilon_m},
\qquad
\varepsilon_i\in\{0,1\},
\]
are precisely the minimal-length representatives of the cosets in $W/S_m$, and every displayed product is length-additive.
\end{lemma}

\begin{proof}
Write a signed permutation as $w=[w_1,\ldots,w_m]$, meaning that $w(e_i)=\operatorname{sgn}(w_i)e_{|w_i|}$ on the standard basis. The reflection $s_i$ interchanges positions $i$ and $i+1$ for $i<m$, and $s_m$ changes the sign in position $m$. In the order
\[
1<2<\cdots<m<-m<\cdots<-1,
\]
each coset modulo $S_m$ contains exactly one increasing signed permutation. It has no right descent among $s_1,\ldots,s_{m-1}$ and is therefore its coset's minimal-length representative.

Right multiplication by $s_j s_{j+1}\cdots s_m$ removes the entry in position $j$ and appends its negative at the end. Starting with the identity, the factors $v_1,\ldots,v_m$ consider the positions $m,m-1,\ldots,1$ in that order. At the step corresponding to $j$, all earlier operations affected only positions greater than $j$, so the entry in position $j$ is still $j$. Selecting that factor changes the sign of $j$; the resulting sequence remains increasing. Thus the displayed products give exactly one increasing representative for each subset of negative values.

If $I$ is the set of negative values in such a representative, its length is
\[
\ell(w)=\sum_{j\in I}(m-j+1).
\]
Indeed, the positive roots are $e_i$ and $e_i\pm e_j$ for $i<j$. No root $e_i-e_j$ becomes negative for an increasing representative. A negative value $j$ contributes one negative image of a root $e_i$, and one negative image of a sum root for each value greater than $j$. This gives the formula. For the displayed product, $I=\{m-i+1:\varepsilon_i=1\}$, so its length is $\sum_i i\varepsilon_i$, the sum of the lengths of its selected factors. The product is therefore length-additive.
\end{proof}

\begin{proposition}\label{prop:binary-bruhat}
Let $G$ be a finite group of Lie type whose relative Weyl group is $W(B_m)$, and let $P$ be the maximal parabolic whose Weyl subgroup is $S_m$. For $v_i$ as in Lemma~\ref{lem:Bweyl}, define
\[
X_i=\{1\}\cup U_{v_i}\dot v_i.
\]
Then multiplication gives an exact product
\[
G=X_1\cdots X_mP.
\]
Each $X_i$ has property $\Dcal_3$.
\end{proposition}

\begin{proof}
Fix a choice of the $\varepsilon_i$ and let $w=v_1^{\varepsilon_1}\cdots v_m^{\varepsilon_m}$. By Lemma~\ref{lem:Bweyl}, this is a reduced product and $w\in W^J$. Successive application of the exact unipotent-product identity shows that the selected sets $U_{v_i}\dot v_i$ multiply bijectively onto a set of the form $U_w\dot wt$, where $t\in T\leq P$. This set contains exactly one representative of every coset in $B\dot wP/P$. Different choices of the $\varepsilon_i$ give different $B$-orbits on $G/P$, not merely different parameters in the same orbit. The disjoint parabolic Bruhat decomposition therefore proves that $X_1\cdots X_m$ is an exact transversal, and hence that $G=X_1\cdots X_mP$ is exact.

Right translation by $\dot v_i^{-1}$ sends $X_i$ to
\[
U_{v_i}\cup\{\dot v_i^{-1}\}.
\]
The element $\dot v_i^{-1}$ does not belong to $U_{v_i}$, since its Weyl-group image is nontrivial. The subgroup $U_{v_i}$ is solvable. Lemma~\ref{lem:solvable-plus-one} and translation invariance give property $\Dcal_3$ for $X_i$.
\end{proof}

\Needspace{12\baselineskip}
\section{Unitary and symplectic groups}

\begin{proposition}\label{prop:unitary}
Every nonabelian finite simple group $G=\PSU_n(q)$ has property $\Dcal_3$.
\end{proposition}

\begin{proof}
Write $n=2m$ or $2m+1$, and let $P$ be the stabilizer of a maximal totally isotropic $m$-space. The relative Weyl group is $W(B_m)$ and the Weyl subgroup of $P$ is $S_m$. Proposition~\ref{prop:binary-bruhat} gives an exact product
\[
G=X_1\cdots X_mP,
\]
with each $X_i$ having property $\Dcal_3$.

Let $Q$ be the unipotent radical of $P$ and $L=P/Q$. The group $Q$ has prime-power order. If $m\geq2$, the Levi quotient $L$ has a normal subgroup $L_0$ which is a central quotient of $\SL_m(q^2)$, and $L/L_0$ is cyclic. Indeed, in even dimension the paired block corresponding to $A\in\GL_m(q^2)$ is $(A^{q\mathsf T})^{-1}$ and the determinant condition is $\det A\in\F_q^{\times}$; in odd dimension the residual one-dimensional block is $(\det A)^{q-1}$, so every $A\in\GL_m(q^2)$ is permitted. Factoring by scalar matrices does not change the conclusion. Proposition~\ref{prop:linear} and Lemma~\ref{lem:normal-extension} give property $\Dcal_2$ for $L_0$, then for $L$, and finally for $P$. When $m=1$, the Levi quotient is solvable and the same conclusion holds. Lemma~\ref{lem:concat} now gives property $\Dcal_3$ for $G$.
\end{proof}

\begin{proposition}\label{prop:symplectic}
Every nonabelian finite simple group $G=\PSp_{2m}(q)$ has property $\Dcal_3$.
\end{proposition}

\begin{proof}
If $m=1$, then $G=\PSp_2(q)\cong\PSL_2(q)$, and the result follows from Proposition~\ref{prop:linear}. Thus assume $m\geq2$. Let $P$ be the stabilizer of a maximal totally isotropic $m$-space. The relative Weyl group is $W(C_m)\cong W(B_m)$ and the Weyl subgroup of $P$ is $S_m$. Proposition~\ref{prop:binary-bruhat} gives
\[
G=X_1\cdots X_mP
\]
exactly, with each $X_i$ having property $\Dcal_3$.

If $Q$ is the unipotent radical of $P$, then $P/Q$ is a central quotient of $\GL_m(q)$. The image of $\SL_m(q)$ is normal with cyclic quotient. Proposition~\ref{prop:linear} and Lemma~\ref{lem:normal-extension} therefore give property $\Dcal_2$ for $P/Q$, and then for $P$. Lemma~\ref{lem:concat} gives property $\Dcal_3$ for $G$.
\end{proof}

\begin{remark}\label{rem:parameters}
Unequal root parameters in the unitary groups do not affect the proof. The binary factors have sizes
\[
q^{2i-1}+1\quad\text{for }\PSU_{2m}(q),
\qquad
q^{2i+1}+1\quad\text{for }\PSU_{2m+1}(q)
\]
for $1\leq i\leq m$, while for $\PSp_{2m}(q)$ they have size $q^i+1$. In odd-dimensional unitary groups, the end subgroup $U_{s_m}$ has order $q^3$ and need not be abelian. Only solvability of the groups $U_{v_i}$ is used.
\end{remark}

\section{Orthogonal groups}

\begin{proposition}\label{prop:odd-orthogonal}
Every nonabelian finite simple group
\[
G=\POmega_{2m+1}(q),\qquad q\text{ odd},
\]
has property $\Dcal_3$.
\end{proposition}

\begin{proof}
If $m=1$, then $G=\POmega_3(q)\cong\PSL_2(q)$, and the result follows from Proposition~\ref{prop:linear}. Thus assume $m\geq2$. Let $P$ stabilize a maximal totally singular $m$-space. The relative Weyl group is $W(B_m)$ and the Weyl subgroup of $P$ is $S_m$. Proposition~\ref{prop:binary-bruhat} gives
\[
G=X_1\cdots X_mP
\]
exactly, with each $X_i$ having property $\Dcal_3$.

Let $Q$ be the unipotent radical of $P$ and $L=P/Q$. On a Witt decomposition, the Levi action contains the block $A\in\GL_m(q)$ together with $A^{-\mathsf T}$ on the paired space. Intersection with $\Omega_{2m+1}(q)$ imposes a spinor-norm condition, but still contains $\SL_m(q)$, and the quotient by this subgroup is abelian. After the projective central quotient, the image of $\SL_m(q)$ remains normal with solvable quotient. Proposition~\ref{prop:linear} and Lemma~\ref{lem:normal-extension} give property $\Dcal_2$ for $L$, then for $P$. Lemma~\ref{lem:concat} gives property $\Dcal_3$ for $G$.
\end{proof}

\begin{proposition}\label{prop:minus-orthogonal}
Every nonabelian finite simple group $G=\POmega^-_{2m}(q)$ has property $\Dcal_3$.
\end{proposition}

\begin{proof}
If $m=2$, then $G=\POmega^-_4(q)\cong\PSL_2(q^2)$, and the result follows from Proposition~\ref{prop:linear}. Thus assume $m\geq3$. Let $r=m-1$ and let $P$ stabilize a maximal totally singular $r$-space. The relative Weyl group is $W(B_r)$ and the Weyl subgroup of $P$ is $S_r$. Proposition~\ref{prop:binary-bruhat} gives
\[
G=X_1\cdots X_rP
\]
exactly, with each $X_i$ having property $\Dcal_3$.

Let $Q$ be the unipotent radical of $P$ and $L=P/Q$. The residual nondegenerate space is an anisotropic plane. In the full orthogonal stabilizer a Levi subgroup is
\[
\GL_r(q)\times O_2^-(q),
\]
where $O_2^-(q)$ is dihedral of order $2(q+1)$ and hence solvable. The intersection defining $\Omega^-_{2m}(q)$ contains $\SL_r(q)\times1$. In odd characteristic the spinor norm on the embedded $\GL_r(q)$ block is the square class of its determinant. In even characteristic, the element $\operatorname{diag}(A,I_2,A^{-\mathsf T})$ has zero Dickson invariant for $A\in\SL_r(q)$. Thus the quotient by $\SL_r(q)\times1$ is solvable. The projective quotient only removes a central subgroup, so the image of $\SL_r(q)$ remains normal with solvable quotient. Proposition~\ref{prop:linear} and Lemma~\ref{lem:normal-extension} give property $\Dcal_2$ for $L$ and then for $P$. Lemma~\ref{lem:concat} gives property $\Dcal_3$ for $G$.
\end{proof}

\begin{remark}
For $\POmega_{2m+1}(q)$ the binary factors have size $q^i+1$ for $1\leq i\leq m$. For $\POmega^-_{2m}(q)$ they have size
\[
q^{i+1}+1\qquad(1\leq i\leq m-1).
\]
Only solvability of the corresponding unipotent parameter groups is needed.
\end{remark}

\begin{proposition}\label{prop:plus}
Every nonabelian finite simple group $G=\POmega^+_{2m}(q)$ has property $\Dcal_6$.
\end{proposition}

\begin{proof}
The case $m=3$ is $\POmega^+_6(q)\cong\PSL_4(q)$ and follows from Proposition~\ref{prop:linear}; smaller dimensions do not define a nonabelian simple plus-type family. Thus assume $m\geq4$ and define
\[
\widehat G=
\begin{cases}
G,& q\text{ even},\\
\PSO^+_{2m}(q),& q\text{ odd}.
\end{cases}
\]
By \cite[Proposition~2.9(ii),(iii)]{BurnessLi}, there is a factorization
\[
\widehat G=HK,
\]
where $H$ is solvable and $K$ is the stabilizer of a nonsingular $1$-space. When $q$ is odd, $G$ has index at most $2$ in $\widehat G$.

Let $J=H\cap K$. Lemma~\ref{lem:solvable-transversal} gives a left transversal $T$ for $J$ in $H$ with property $\Ecal$. Hence
\[
\widehat G=TK
\]
is exact: $H=TJ$ and $J\leq K$, while $t_1k_1=t_2k_2$ with $t_i\in T$ and $k_i\in K$ implies $t_2^{-1}t_1\in H\cap K=J$, so $t_1=t_2$. If $q$ is even, every nonabelian composition factor of $K$ belongs to the symplectic family of rank $m-1$ and has property $\Dcal_3$ by Proposition~\ref{prop:symplectic}. If $q$ is odd, every nonabelian composition factor of $K$ belongs to the odd-dimensional orthogonal family and has property $\Dcal_3$ by Proposition~\ref{prop:odd-orthogonal}. All remaining composition factors are cyclic. Repeated use of Lemma~\ref{lem:normal-extension} gives property $\Dcal_3$ for $K$. Since $T$ has property $\Dcal_2$, Lemma~\ref{lem:concat} gives property $\Dcal_3$ for $\widehat G$.

If $q$ is even, $\widehat G=G$. If $q$ is odd, Lemma~\ref{lem:index} transfers the property to the index-at-most-two subgroup $G$, with constant at most $6$.
\end{proof}

In even characteristic, $\Omega_{2m+1}(q)\cong\operatorname{Sp}_{2m}(q)$, so Proposition~\ref{prop:symplectic} covers the odd-dimensional family there as well. Thus Propositions~\ref{prop:odd-orthogonal}, \ref{prop:minus-orthogonal}, and~\ref{prop:plus} cover all orthogonal simple groups.

\section{Exceptional groups}

The exceptional families have bounded relative rank, so a much coarser Bruhat product is sufficient.

\begin{proposition}\label{prop:exceptional}
There is an absolute constant $C_{\rm exc}$ such that every finite simple exceptional group of Lie type has property $\Dcal_{C_{\rm exc}}$.
\end{proposition}

\begin{proof}
For a member of an exceptional family carrying its standard split or twisted $BN$-pair, let $(B,\mathcal N)$ be the standard $BN$-pair with relative Weyl group $W$, and write $\ell$ for the Coxeter length function on $W$. Let $w_0$ be the longest element, set
\[
\ell_0=\ell(w_0),
\]
and choose a reduced expression
\[
w_0=s_1\cdots s_{\ell_0}.
\]
For $1\leq i\leq\ell_0$, let $P_i$ be the minimal parabolic corresponding to $s_i$. Let $U_i=U_{s_i}$ be the associated unipotent subgroup, choose a representative $\dot s_i\in\mathcal N$, and define
\[
X_i=\{1\}\cup U_i\dot s_i.
\]
Then $X_i$ represents the left cosets of $B$ in $P_i$, so $P_i=X_iB$ is exact. Here $p$ is the defining characteristic, $U_i$ is a $p$-group, and
\[
|X_i|=q_i+1,
\qquad q_i=|U_i|\geq2.
\]
Every element of $W$ is below $w_0$ in Bruhat order, and the subword property shows that it is represented by a subword of the fixed reduced expression for $w_0$. Thus every Bruhat cell is contained in the corresponding product of minimal parabolics, and
\[
G=P_1\cdots P_{\ell_0}.
\]
Since $P_i=X_iB$ and $BX_iB=P_i=X_iB$, the intervening copies of $B$ may be absorbed successively to the right. Hence
\[
G=X_1\cdots X_{\ell_0}B.
\]
Each $X_i$ has property $\Dcal_3$ by Lemma~\ref{lem:solvable-plus-one}; the Borel subgroup $B$ is solvable and has property $\Dcal_2$.

Let
\[
M=|B|\prod_{i=1}^{\ell_0}(q_i+1).
\]
The longest Bruhat cell has size
\[
|B|\prod_{i=1}^{\ell_0}q_i.
\]
The multiplication map from $X_1\times\cdots\times X_{\ell_0}\times B$ onto $G$ is surjective, so $M\geq|G|$. Hence, with $\delta=M/|G|$,
\[
1\leq\delta\leq
\prod_{i=1}^{\ell_0}\left(1+\frac1{q_i}\right)
\leq\left(\frac32\right)^{\ell_0}.
\]
For exceptional groups one has $\ell_0\leq120$, with equality only for type $E_8$. For example, the longest-element lengths are $6,24,36,63,120$ for $G_2,F_4,E_6,E_7,E_8$ respectively; the twisted exceptional relative Weyl groups have no larger length.

Let $x_1\cdots x_r=|G|$. The case $r=1$ is immediate. For $r\geq2$, apply Lemma~\ref{lem:concat} to the product $X_1\cdots X_{\ell_0}B$, whose size parameter is $M$, with the numbers
\[
x_1,\ \delta x_2,\ x_3,\ldots,x_r.
\]
These numbers have product $M$. Lemma~\ref{lem:concat} gives a covering of $G$ in which the first subset has size at most $x_1$, the second has size at most $3\delta x_2$, and every later subset has size at most $3x_i$. Hence every group covered by the standard $BN$-pair argument has property $\Dcal_{C_{\rm exc}}$ with, for example,
\[
C_{\rm exc}=3\left(\frac32\right)^{120}.
\]
Only finitely many simple exceptional groups arise from the exceptional small parameters where the standard ambient group is not itself simple. Each fixed finite group has property $\Dcal_C$ for some $C$, so enlarging $C_{\rm exc}$ once covers this finite list as well. In particular, the Tits group may alternatively be obtained from ${}^2F_4(2)$ by the bounded-index transfer in Lemma~\ref{lem:index}.
\end{proof}

This includes the Suzuki, small Ree, large Ree, triality, and all untwisted exceptional families. No separate treatment of the first few field sizes is needed.

\section{Proof of the main theorem}

\begin{proof}[Proof of Theorem~\ref{thm:lie}]
Proposition~\ref{prop:linear} covers the linear family. Propositions~\ref{prop:unitary} and~\ref{prop:symplectic} cover the unitary and symplectic families. Propositions~\ref{prop:odd-orthogonal}, \ref{prop:minus-orthogonal}, and~\ref{prop:plus} cover the orthogonal families, including the even-characteristic odd-dimensional case through the standard symplectic isomorphism. Proposition~\ref{prop:exceptional} covers the exceptional families. The usual low-rank isomorphisms introduce no additional infinite family, and the finitely many small exceptional parameters are absorbed into the constant in Proposition~\ref{prop:exceptional}. Thus every nonabelian finite simple group of Lie type has property $\Dcal_{C_{\rm Lie}}$ for
\[
C_{\rm Lie}=\max\{6,C_{\rm exc}\}.
\]
For prescribed weights $\alpha_1,\ldots,\alpha_h$, take $x_i=|G|^{\alpha_i}$ in Definition~\ref{def:C}. This gives the conclusion of Theorem~\ref{thm:lie}.
\end{proof}

\begin{proof}[Proof of Theorem~\ref{thm:global}]
By the classification of finite simple groups, a finite simple group is cyclic of prime order, alternating, of Lie type, or sporadic; see, for example, \cite{Wilson}. Cyclic simple groups are solvable and have property $\Dcal_2$. Alternating simple groups have property $\Dcal_2$ by \cite[Theorem~1]{KozmaLev}. The finite simple groups of Lie type have property $\Dcal_{C_{\rm Lie}}$ by Theorem~\ref{thm:lie}. Since there are only finitely many sporadic simple groups, one constant works for all of them: for a fixed sporadic simple group $S$, property $\Dcal_{|S|}$ follows by taking the first subset to be $\{1\}$, the second to be $S$, and all later subsets to be $\{1\}$ whenever at least two factors occur. Hence there is an absolute constant $C_0$ such that every finite simple group has property $\Dcal_{C_0}$.

Induction on $|G|$ now gives the result for every finite group. The trivial group is immediate. If $G$ is nontrivial and simple, the preceding paragraph applies. Otherwise choose a nontrivial proper normal subgroup $N\lhd G$. By induction, both $N$ and $G/N$ have property $\Dcal_{C_0}$, and Lemma~\ref{lem:normal-extension} gives property $\Dcal_{C_0}$ for $G$. Given $\alpha_1,\ldots,\alpha_h$ as in the theorem, take $x_i=|G|^{\alpha_i}$ in Definition~\ref{def:C}. Since the covering subsets lie in $G$, their product is exactly $G$.
\end{proof}

\begin{proof}[Proof of Theorem~\ref{thm:rohrbach}]
Apply Theorem~\ref{thm:global} with $\alpha_1=\cdots=\alpha_h=1/h$ and let
\[
A=A_1\cup\cdots\cup A_h.
\]
Then $A^h=G$ and
\[
|A|\leq \bigl(1+C(h-1)\bigr)|G|^{1/h}.
\]
Thus one may take $c_h=1+C(h-1)$.
\end{proof}

\section*{Acknowledgements}
The author thanks his colleague Xingde Jia, who introduced him to this problem about fifteen years ago. The author acknowledges support from the Simons Foundation Award \#918096. OpenAI's ChatGPT was used for mathematical exploration and editorial assistance. The author checked the final arguments and is responsible for the contents of the paper.

\end{document}